\documentclass[12pt, a4paper]{article}
\usepackage[T1]{fontenc}
\usepackage{graphicx} % Required for inserting images
\usepackage{amsmath}
\usepackage{amsthm}
\usepackage{xcolor}
\title{On the lower- and upper bounds of three-colored anti-Ramsey numbers}
\author{Máté Giesz\\mgiesz17000@gmail.com}
\date{2026}
\newtheorem{defi}{Definition}[section]
\newtheorem{prop}{Proposition}[section]
\newtheorem{thm}{Theorem}[section]
\newtheorem{cor}{Corollary}[section]
\newtheorem{rem}{Remark}[section]
\usepackage{url}
\usepackage{amssymb}
\usepackage{setspace}
\usepackage{mathrsfs}
\usepackage{mathabx}
\usepackage{amsfonts}
\usepackage[colorlinks,pagebackref]{hyperref}
\hypersetup{
    linkcolor=black!70!black,  % A tartalomjegyzék és a belső hivatkozások (képletek, ábrák) színe
    citecolor=green,  % Az irodalomjegyzékre mutató hivatkozások színe [1]
    urlcolor=blue     % A külső webes linkek színe
}
\renewcommand*{\backref}[1]{}
\renewcommand*{\backrefalt}[4]{%
  \ifcase #1 % Ha nincs hivatkozás
    (Not cited.)%
  \or % Ha egyetlen hivatkozás van
    $\uparrow$#2. p.%
  \else % Ha több hivatkozás van
    $\uparrow$#2. pp.%
  \fi
}

\begin{document}

\maketitle
\newpage
\tableofcontents
\newpage
\begin{abstract}
    In this paper we study a relatively new Ramsey type graph parameter introduced by my Bachelor's thesis supervisor Nóra Almási, and her Master's thesis supervisor Maria Axenovich. \cite{almasi2025online} The parameter $\xi(G)$ denotes the largest natural number $n$, such that there is a red-blue-green edge coloring of $K_n$, such that in the edge coloring all copies of $G$ represent all three colors. After a brief recap on Ramsey numbers, we precisely define $\xi(G)$, and take a look at its relations to the Ramsey numbers. In the third chapter we establish lower bounds on $\xi(G)$, in the cases $G=K_5$ and $G=K_6$, namely $\xi(K_5)\ge 19$ and $\xi(K_6)\ge 31$.
\end{abstract}
\section{Introduction}
Ramsey theory is a remarkable part of combinatorics, especially graph theory. This research field has several open questions, despite being almost $100$ years old.\\
In its simplest form, Ramsey's theorem states the following. For every pair of simple, undirected graphs $G$ and $H$ there exists a least integer $R=R(G,H)$, such that every red-blue edge coloring of the complete graph $K_R$ on $R$ vertices has either a red copy of $G$, or a blue copy of $H$.
\begin{defi}
    We call such least $R=R(G,H)$ the \emph{two-color Ramsey number corresponding to $G$ and $H$}. When $G=K_k$ and $H=K_l$, we use the notation $R(k,l)$ for brevity.
\end{defi}
\begin{prop}
    One can easily prove that $R(G,H)\le R(n_1,n_2)$, where $n_1=|V(G)|$, and $n_2=|V(H)|$.
\end{prop}
These definitions seem very obscure at first glance. One may ask the agonizing question: \textit{"Do all the Ramsey numbers exist?"} By the results from two prominent figures of Ramsey theory, Pál Erdős and György Szekeres \cite{erdos1935combinatorial}, we know the following results concerning the Ramsey numbers $R(k,l)$.
\begin{thm}
    \emph{\cite{erdos1935combinatorial}} For the Ramsey number $R(k,l)$ the following upper bound holds: $$R(k,l)\le R(k-1,l)+R(k,l-1),$$ from which one can prove $$R(k,l)\le\binom{k+l-2}{l-1}=\binom{k+l-2}{k-1}$$ by induction on $k$ and $l$. When $k=l$, it follows that $$R(k,k)\le4^{k-1},$$ from the definition of binomial coefficients.
\end{thm}
Note that the last upper bound is not sharp. In $2024$ the group of P.~Gupta, N.~Ndiaye, S.~Norin, and L.~Wei improved the upper bound of $R(k,k)$ to $3.8^{k+o(k)}$. \cite{gupta2024optimizing}\\
The following lower bound was proved by Erdős in $1947$. \cite{erdos1947some} Similarly to the previous upper bounds, this result was important, since it introduced a new type of proof, named the \textit{probability method}.
\begin{thm}\emph{(Erdős, $1947$)} 
    For sufficiently large $k\in\mathbb{N}$, the lower estimate $R(k,k)>\left(\sqrt{2}\right)^k$ holds.
\end{thm}
\begin{cor}
    $R(k,l)$ is well defined, as from Theorem 1.1, $R(k,l)\le\binom{k+l-2}{k-1}<\infty.$
    By Proposition 1.1 $R(G,H)<\infty$ also holds, making these types of Ramsey numbers well defined as well.
\end{cor}
\section{The exact definition and well-definedness of the three-colored anti-Ramsey number $\xi(G)$}
The paper will inspect the following graph parameter denoted by $\xi(G)$. This parameter is closely related to the \textit{Online Ramsey Turnaround Game}. This game was mainly studied by Nóra Almási and her Master's Thesis supervisor, Maria Axenovich \cite{almasi2025online}, as well as me in my Bachelor's Thesis. \cite{giesz2026ramsey} Regardless, we will look at this parameter in purely "Ramsey theory context", meaning we do not deal with the game.
\begin{defi}
    \emph{(three-colored anti-Ramsey number)} Let $G$ be a simple, undirected graph. Then $\xi(G)$ denotes the biggest positive integer $n$, with the following property. There exists a red-blue-green edge-coloring of $K_n$, such that in every $G\subseteq K_n$ copy of $G$, all three colors get represented. (More precisely in such edge-coloring, every subgraph of $K_n$ which is isomorphic to $G$ has at least one red edge, and one blue edge, and one green edge.) Similarly to Ramsey number $R(k,l)$, we denote $\xi(K_k)$ by $\xi(k)$.
\end{defi}
Similarly to the two-color Ramsey numbers, one may ask the agonizing question: \textit{"Does $\xi(G)$ exist for every graph $G$?"} By the following theorem concerning $\xi(k)$, we will see that the answer is \textit{yes} to the previous question.
\begin{thm}
    For every $k\ge 3$, $\xi(k)<R(k,k)$ holds.
\end{thm}
\begin{proof}
    We need to show that in every red-blue-green edge-coloring of $K_{R(k,k)}$ there is a copy of $K_k$ representing at most two colors. Let $n=R(k,k)$, and consider an arbitrary red-blue-green edge-coloring of $K_n$. We call this edge-coloring the \textit{original} coloring. From this we define an auxiliary \textit{second} edge-coloring on $K_n$. Let every red and blue edge be \textit{yellow}, and let every green edge remain green in the second coloring. This way we can apply Ramsey's theorem on $K_n$, meaning there is either a yellow or green copy of $K_k$. (More precisely a subgraph isomorphic to $K_k$, with all edges being yellow, or all edges being green.) In the original coloring, the yellow $K_k$ means a $K_k$ subgraph, having only red and blue edges, while the green $K_k$ remains a green $K_k$. This proves the theorem.
\end{proof}
\begin{cor}
    One can easily verify that $\xi(G)\le\xi(n)$, where $n=|V(G)|$. Therefore the previous theorem yields that $\xi(G)$ is also well defined.
\end{cor}
The following proposition is an easy observation, that yields an improvement for the upper bound of $\xi(k)$. \cite{almasi2026personal} We leave out the proof of this proposition, as it is similar to the proof of Theorem 2.1.
\begin{prop}
    For every $k\ge 3$, $\xi(k)<R(K_k,K_k-e)$ holds, where $K_k-e$ denotes the simple graph with $k$ vertices and $\binom{k}{2}-1$ edges.
\end{prop}
\begin{cor}
    \emph{(Upper bounds of $\xi(5)$ and $\xi(6)$) }By Proposition 2.1, we  have $\xi(5)\le 29$ and $\xi(6)\le103$.
\end{cor}
\begin{proof}
    In $2026$ V.~Angeltveit showed that $R(K_5,K_5-e)=30$ \cite{angeltveit2026exact}, meaning $\xi(5)<30$, while in $2012$ L. Boza showed that $R(K_6,K_6-e)\le104$ \cite{boza2012numeros}, meaning $\xi(6)<104$. 
\end{proof}
\section{Lower bounds for $\xi(5)$ and $\xi(6)$}
The following previously known theorem states an easily verifiable lower bound on $\xi(n)$, using a classical tool in combinatorics called \textit{pigeonhole principle}. However more importantly the lower bound puts into perspective, how much improvement our constructions exactly mean.
\begin{thm}
    For every $n\ge 3$, $\xi(n)\ge(n-1)^2$ holds, moveover there is equality in the case of $n=3$ and $n=4$.
\end{thm}
\begin{proof}
    Let $y=(n-1)^2$, and let $V(K_y)=\{(a,b)\in\mathbb{Z}^2\mid 1\le a,b\le n-1\}$. Let the color of the edges $(a,b)(a,b'), b\ne b'$ be red, $(a,b)(a',b), a\ne a'$ be blue, and the remaining edges be green. That is, we imagine $V(K_y)$ as an $(n-1)\times(n-1)$ grid, the \textit{"horizontal"} edges are colored red, the \textit{"vertical"} edges are colored blue, and the \textit{"diagonal"} edges are colored green. This way, by the \textit{pigeonhole principle}, every $K_n\subseteq K_y$ has an edge of each type, that is, all three colors get represented in the $K_n$ copies.\\
    We prove the second part only when $n$ is equal to $3$. The result for $n=4$ was proved by András Gyárfás, as well as K.M.~Chung and C.L.~Liu. \cite{chung1978generalization} By the first part $\xi(3)\ge(3-1)^2=2^2=4$. Consider an arbitrary $3$-edge-coloring of $K_5$, and let $v$ be an arbitrary vertex of the graph. By the \textit{pigeonhole principle} at least two edges sharing the common vertex $v$, say $va$ and $vb$ must share the same color. Therefore $vab$ is a $K_3$ copy representing at most two colors, yielding $\xi(3)<5$.
\end{proof}
\begin{cor}
    Substituting $n=5$ and $n=6$ in the previous inequality, we get $\xi(5)\ge16$ and $\xi(6)\ge 25$.
\end{cor}
The following construction was established by Christian Ortlieb in his Bachelor's Thesis, in $2017$. \cite{ortlieb2017balanced} He stated and proved a lemma characterizing $3$-edge-colorings of $K_n$ containing a $K_k$ subgraph representing at most two colors. He checked the conditions of this lemma using Matlab, meaning that the lower bound had no elementary proof or elementary way of being discovered. In the proof of the following theorem we construct the ominous $3$-edge-coloring, and analytically show that the construction is indeed rigorous.
\begin{thm}
    The $3$-edge-coloring of $K_{19}$ defined below fulfills our conditions, namely every copy of $K_5\subseteq K_{19}$ represents all three colors, implying $\xi(5)\ge 19$.
\end{thm}
\begin{proof}
    Let the vertices of $K_{19}$ be the numbers $0,1,\dots,18$. As $19=6k+1$, where $k=3$, we will compute the \textit{cubic residues} modulo $19$. If we compute the cubes of non-zero elements, one can see that the set of cubes decomposes into $6$ classes. Let
    \begin{align*}
        H &= \{1, 7, 8, 11, 12, 18\} \equiv \{\pm1, \pm7, \pm8\}, \\
        2H &= \{2, 3, 5, 14, 16, 17\} \equiv \{\pm2, \pm3, \pm5\}, \text{ and} \\
        4H &= \{4, 6, 9, 10, 13, 15\} \equiv \{\pm4, \pm6, \pm9\}.
    \end{align*}
    (By $kH$ we mean the elementwise multiplication of $H$ with the number $k$, modulo $19$. Here $k=2$ and $k=4.$) Let the color of the edge $xy, 0\le x<y\le18$ be red if $x-y\in H$, green if $x-y\in 2H$, and blue if $x-y\in 4H$. Notice that as $-1$ is a cubic residue, if $x-y\in H$, then $y-x\in H$ meaning the colors of the edges are undirected. It is also clear from the definition, that the sets $H,2H,$ and $4H$ form a partition of the set $[18]$, meaning that every edge of the $K_{19}$ indeed gets a color.\\
    Next we prepare the proof by introducing some notations simplifying the structure of the proof. Let $G_r,G_b,G_g$ denote the subgraphs of $K_{19}$ determined by the red, blue, and green edges respectively. Notice that a $K_5$ subgraph containing at most two colors in $K_{19}$ means that either $G_r,G_b,G_g$ contains five independent vertices, or in other words the \textit{clique number} of the complements of these graphs is at least $5$. Another important observation is that the colors play a symmetric role, therefore it is enough to prove that $\alpha(G_r)=\omega(\overline{G_r})\le 4$.\\
    Assume for contradiction that $K_5\subseteq\overline{G_r}$. Organize the vertices of $K_{19}$ in a circle, and let the distances of the $K_5$ subgraph be denoted by $d_1,d_2,\dots,d_5$. This way we obtain $\sum\limits_{i=1}^5d_i=19$. As the forbidden distances in $\overline{G_r}$ are $1,7,8$ and their negatives modulo $19$, we have $d_i\ge 2$ and $d_i\ne7,8$. We also notice that any proper subset of vertices of the $K_5$ also forms a clique, thus the sum of neighbouring terms in the sum above cannot be equal to $1,7,8$ either. This yields that up to permutation, the following tuples cannot be neighbouring terms in the above partitions of $19$: $(2,6),(3,5),(4,4),(2,5),(3,4),(2,2,3)$ $(2,2,4),(2,3,3),(2,2,2,2),(2,4,2,2)$.\\
    Next we consider the largest term in the partitions of $19$, without loss of generality we can assume that this term is $d_1$.\\
    If $d_1\ge 13$ we have $8\le d_2+\cdots+d_5\le 6$, which is a contradiction. It is also clear that $d_1=11,12$ are also impossible, as $11$ and $12$ are the negatives of $8$ and $7$ modulo $19$.\\
    If $d_1=10$, we get $d_2+\cdots+d_5=9$ which yields that the only possible partition of $19$ is $19=10+2+2+2+3$, a contradiction, as the triplet $(2,2,3)$ is forbidden up to permutation.\\
    The case $d_1=9$ yields that $d_2+\dots+d_5=10$. This sum cannot contain any term greater than $4$, as this would mean that the sum of the remaining $3$ terms would be at most $5$. Therefore the only possible partitions are $19=9+2+3+2+3=9+2+2+2+4$ and their permutations. Both of these configurations are forbidden, because of the triplets $(2,2,4),(2,2,3)$ and $(2,3,3)$.\\
    The next case is $d_1=6$, as $d_i=7,8$ is forbidden. In this case we have $d_2+\cdots+d_5=13$. If the partition had another term equal to $5$ or $6$, the remaining three numbers could not be equal to $2$, since $5+2=7$ and $6+2=8$ are forbidden. This yields that $\sum\limits_{i=1}^5d_i\ge 5+4+4+4+3=20$, a contradiction. Therefore the only possible partitions are $19=6+4+3+3+3=6+4+3+4+2$, which are both forbidden because of the pairs $(4,4),(4,2)$ and $(4,3)$.\\
    Next we inspect $d_1=5$. This means that $d_2+\cdots+d_5=14$ The partitions of $19$ having another term equal to $5$ would imply that the remaining three terms cannot be equal to $2$ or $3$, since $5+2=7$ and $5+3=8$. This implies that $\sum\limits_{i=1}^5\ge5+4+4+4+4=21$, which is not possible. Therefore the only possible partitions are $19=5+4+4+4+2=5+4+4+3+3$. Neither of them is possible, as $(4,4)$ and $(4,3)$ are forbidden.\\
    Finally we inspect some trivial cases. If $d_1=4$, we obtain $d_2+\cdots+d_5=15$, which is only possible if $19=4+4+4+4+3$, contradicting with $(4,3)$ being forbidden. The case $d_1\le3$ yields $\sum\limits_{i=1}^5d_i\le15<19$, which is a contradiction.\\
    Overall our assumption $K_5\subseteq\overline{G_r}$ was false, which proves the theorem.
\end{proof}
In the proof of the next theorem we present an analogous construction on $K_{31}$, concerning the lower bound for $\xi(6)$. This construction had no previous proofs, regarding if it fills the given criteria of $\xi(6)$.
\begin{thm}
    The $3$-edge coloring of $K_{31}$ defined below fulfills our conditions, namely every copy of $K_6\subseteq K_{31}$ represents all three colors, implying $\xi(6)\ge 31$.
\end{thm}
\begin{rem}
    The proof of Theorem 3.2 follows the same structure as the previous proof, as well as the construction of the edge coloring of $K_{31}$ compared to $K_{19}$. A potential research topic for the future could be if these type of constructions can be generalized for bigger complete graphs, improving the general lower bound of $\xi(n)\ge(n-1)^2$.
\end{rem}
\begin{proof}
    (Theorem 3.3) 
    Let the vertices of $K_{31}$ be the numbers $0, 1, \dots, 30$. As $31=6k+1$, where $k=5$, we will compute the \textit{cubic residues} modulo $31$. If we compute the cubes of non-zero elements, one can see that the set of cubes decomposes into $10$ classes. Let 
    \begin{align*}
        H &= \{1, 2, 4, 8, 15, 16, 23, 27, 29, 30\} \equiv \{\pm1, \pm2, \pm4, \pm8, \pm15\}, \\
        3H &= \{3, 6, 7, 12, 14, 17, 19, 24, 25, 28\} \equiv \{\pm3, \pm6, \pm7, \pm12, \pm14\}, \text{ and} \\
        9H &= \{5, 9, 10, 11, 13, 18, 20, 21, 22, 26\} = \{\pm5, \pm9, \pm10, \pm11, \pm13\}.
    \end{align*}
    (By $kH$ we mean the elementwise multiplication of $H$ with the number $k$, modulo $31$. Here $k=3$ and $k=9$.) Let the color of the edge $xy, 0\le x<y\le30$ be red if $x-y\in H$, green if $x-y\in 3H$, and blue if $x-y\in 9H$. Notice that as $-1$ is a cubic residue, if $x-y\in H$, then $y-x\in H$ meaning the colors of the edges are undirected. It is also clear from the definition, that the sets $H, 3H,$ and $9H$ form a partition of the set $[30]$, meaning that every edge of the $K_{31}$ indeed gets a color.\\
    For the upcoming parts of the proof we keep the notations used in the previous proof. This way we need to show that $\omega(\overline{G_r})\le 5$. If we assume by contradiction that $K_6\subseteq\overline{G_r}$, by a similar argument to the previous proof, one can see that $\sum\limits_{i=1}^6d_i=31$. Since $d_i=1,2$ is forbidden, we obtain $d_i\ge 3$. Again  consider the largest term in the above partitions of $31$, without loss of generality, we can assume that this term is $d_1$.\\
    The case $d_1\ge 17$ yields $15\le d_2+\cdots+d_6\le 14$, a contradiction. The case $d_1=16$ cannot hold, as $16=-1\cdot 15$ modulo $31$, the case $d_1=15$ cannot hold by definition.\\
    If $d_1=14$, then $d_2+\cdots+d_6=17$. If this sum had a term at least $6$, the remaining four numbers would add up to at most $11$, which is impossible. This yields that the only possible partitions are $31=14+4+4+3+3+3=14+5+3+3+3+3$. Neither of them are possible, as $4$ and $5+3=8$ are forbidden.\\
    If $d_1=13$, we have $d_2+\cdots+d_6=18$. If this sum had a term with value at least $7$, we would have that the sum of the remaining four numbers is at most $11$, a contradiction. If the partition had a term equal to $6$, we would have that the remaining four numbers add up to $12$, yielding the partition $31=13+6+3+3+3+3$, which is not possible, as $6+3+3+3=15$ is forbidden. This implies that the only possible partitions up to permutations are $31=13+5+4+3+3+3=13+4+4+4+3+3$. Both of these partitions yield a contradiction, as the term $4$ is forbidden.\\
    The next case is $d_1=11$, as $d_1=12$ is forbidden. In this case we obtain $d_2+\cdots+d_6=20$. If a term other than $d_1$ was at least $9$, the remaining four numbers would add up to $11$, giving us a contradiction. This yields that the remaining terms are at most $7$, as the term $d_i=8$ is forbidden. If the partition has a term $7$, we obtain the partition $31=11+7+3+3+3+4$, a contradiction, as $4$ is forbidden. Otherwise the possible partitions up to permutation are $31=11+6+5+3+3+3=11+6+4+4+3+3=11+5+5+4+3+3=11+4+4+4+4+4$. Neither of them is possible, as $5+3=8$, $6+3+3+3=15$, and $4$ are forbidden respectively.\\
    If $d_1=10$, then $d_2+\cdots+d_6=21$. If another term excluding $d_1$ had a value of $10$, the remaining four terms would add up to $11$, which yields a contradiction. If a term had a value of $9$, the remaining four terms would have a sum of $12$, yielding the only possibility of $31=10+9+3+3+3+3$, which is a contradiction, as $9+3+3=15$. Otherwise the only suitable partitions are $31=10+7+5+3+3+3=10+6+6+3+3+3=10+6+5+4+3+3=10+6+4+4+4+3=10+5+4+4+4+4$, all five of them leading to a contradiction. The reason behind this is that $4$, $5+3=8$, $6+3+3=12$, and $6+3+3+3=15$ are forbidden.\\
    The case $d_1=9$ yields $d_2+\cdots+d_6=22$. If the partition had another term other than $d_1$ with the value $9$, the remaining four terms would add up to $13$. This implies that the only partition that works is $31=9+9+4+3+3+3$, which does not work, as $4$ is forbidden. Since $d_i\ne8$, we have that the only term of the partition with value $9$ is $d_1$, and the other terms have a value of at most $7$. This means that the only possible partitions up to permutations are $31=9+7+3+3+3+6=9+7+3+3+4+5=9+6+4+3+4+5=9+5+5+4+4+5$, which all lead to a contradiction, as $6+3+3+3=15$, $6+3+3=12$, and $4$ are forbidden.\\
    The distance $d_1=8$ is not allowed, therefore we move onto $d_1=7$, meaning $d_2+\cdots+d_6=24$. If any of the remaining five terms had a value of $7$, the remaining four terms would have a sum of $17$, yielding the partitions $31=7+7+4+4+4+5=7+7+4+4+3+5$, which are forbidden, since $4$ is forbidden. Otherwise the following partitions, and their permutations remain: $31=7+4+5+5+5+5=7+6+3+5+5+5=7+6+3+4+6+5=7+6+4+6+5+3=7+6+4+6+4+4$. Either of them yields a contradiction, since $4$, $5+3=8$, and $6+6+3=15$ are forbidden.\\
    If $d_1=6$, we have $d_2+\cdots+d_6=25$. This implies that the only possible partitions are $31=6+5+5+5+5+5=6+6+5+5+5+4=6+6+5+5+4+4=6+6+6+5+4+4$. As $4$, $5+3=8$, and $4+5+6=5+5+5=15$ are forbidden, these partitions, and their permutations yield a contradiction.
    In the end we inspect the case $d_1\le 5$. In this case we get $\sum\limits_{i=1}^6d_i\le30<31$, a contradiction, as the sum has to be equal to \textit{exactly} $31$.\\
    Again, overall our assumption $K_5\subseteq\overline{G_r}$ was false, which proves the theorem.
\end{proof}
\section{Summary}
In the paper we introduced, and explored the graph parameter $\xi(G)$ in the case of $G=K_5$ and $G=K_6$. After introducing and proving the well known lower bound of the parameter on complete graphs, we made edge-coloring constructions for $K_{19}$ and $K_{31}$, concerning the parameter $\xi(5)$ and $\xi(6)$. The first construction was discovered and proved previously by Christian Ortlieb, however it had no elementary proof or discovery prior to our proof. The second construction was completely unknown, however it can be considered a generalization of the first construction, while both of them are generalizations of the so called \textit{Paley graphs}, mainly studied by Robert Greenwood and Andrew Gleason in $1955$. \cite{greenwood1955combinatorial} A potential future research topic would be if these type of constructions can be generalized for bigger cliques. Finding the exact values of $\xi(n)$ can also be a future research topic. While $\xi(n)$'s value is known in the case $n=3,4$, for $n\ge 5$ only general bounds hold regarding the graph parameter. We conjecture that the latter has a similar difficulty as finding exact values for Ramsey number $R(k,k)$.
\begin{quote}
  \textit{„Suppose aliens invade the Earth and threaten to obliterate it in a year's time unless human beings can find the Ramsey number $R(5,5)$. We could marshal the world's best minds and fastest computers, and within a year we could probably calculate the value. If the aliens demanded the Ramsey number $R(6,6)$, however, we would have no choice but to launch a preemptive attack.”} \\
  \vspace{0.2cm}
  -- Pál Erdős
\end{quote}

\end{document}